\documentclass[12pt]{amsart}

\usepackage{times,epsf}

\begin{document}

\newtheorem{thm}{Theorem}[section]
\newtheorem{lem}[thm]{Lemma}
\newtheorem{cor}[thm]{Corollary}

\theoremstyle{definition}
\newtheorem{defn}{Definition}[section]

\theoremstyle{remark}
\newtheorem{rmk}{Remark}[section]

\def\square{\hfill${\vcenter{\vbox{\hrule height.4pt \hbox{\vrule
width.4pt height7pt \kern7pt \vrule width.4pt} \hrule height.4pt}}}$}

\def\T{\mathcal T}

\newenvironment{pf}{{\it Proof:}\quad}{\square \vskip 12pt}

\title[Minimizing CMC Hypersurfaces in Hyperbolic Space]{Minimizing Constant Mean Curvature Hypersurfaces in Hyperbolic Space}
\author{Baris Coskunuzer}
\address{Department of Mathematics \\ Yale University \\ New Haven, CT 06520}
\email{baris.coskunuzer@yale.edu}

\maketitle


\newcommand{\cirD}{\overset{\circ}{D}}
\newcommand{\Si}{S^2_{\infty}({\Bbb H}^3)}
\newcommand{\SI}{S^n_{\infty}({\Bbb H}^{n+1})}
\newcommand{\PI}{\partial_{\infty}}
\newcommand{\Hyp}{{\Bbb H}^{n+1}}
\newcommand{\BH}{\Bbb H}
\newcommand{\BR}{\Bbb R}
\newcommand{\BC}{\Bbb C}
\newcommand{\BZ}{\Bbb Z}

\begin{abstract}

We study the constant mean curvature (CMC) hypersurfaces in $\Hyp$ whose asymptotic boundaries are closed
codimension-$1$ submanifolds in $\SI$. We consider CMC hypersurfaces as generalizations of minimal
hypersurfaces. We naturally generalize some notions of minimal hypersurfaces like being area minimizing,
convex hull property, exchange roundoff trick to the CMC hypersurface context. We also give a generic
uniqueness result for CMC hypersurfaces in hyperbolic space.

\end{abstract}

\section{Introduction}

In this paper, we are interested in complete CMC hypersurfaces in $\Hyp$ whose asymptotic boundaries are
closed codimension-$1$ submanifolds in $\SI$. CMC hypersurfaces in hyperbolic space became attractive for
researchers after the impressing results on minimal hypersurfaces in hyperbolic space in 1980's. In [A1],
Anderson showed the existence of area minimizing hypersurfaces in $\Hyp$ for any given codimension-$1$ closed
submanifold in $\SI$. He also showed that for mean convex asymptotic boundaries, there exist a unique minimal
hypersurface in hyperbolic space. Then Hardt and Lin studied the regularity of these area minimizing
hypersurfaces in [HL]. They also generalized Anderson's uniqueness result to star-shaped domains in $\SI$.

After these results on minimal hypersurfaces in hyperbolic space, the question of generalization of these
results to CMC hypersurfaces was naturally arised. In the following decade, there have been many important
generalizations of these results to CMC hypersurfaces in hyperbolic space. In [To], Tonegawa generalized
Anderson's existence and Hardt and Lin's regularity results for CMC hypersurfaces by using geometric measure
theory methods. In the same year, by using similar techniques, Alencar and Rosenberg got a similar existence
result in [AR]. By using analytic techniques, Nelli and Spruck generalized Anderson's uniqueness result for
mean convex domains to CMC context in [NS]. Then, Guan and Spruck extended Hardt and Lin's uniqueness results
for star-shaped domains to CMC hypersurfaces in hyperbolic space in [GS].

As above paragraphs suggest it is very useful to think CMC hypersurfaces as generalizations of minimal
hypersurfaces. Many techniques for area minimizing hypersurfaces can be applied to the CMC hypersurfaces. In
this paper, we will focus on this point, and give natural generalizations of some notions about minimal and
area minimizing hypersurfaces to the CMC hypersurfaces.

The main results of the paper are as follows. First, we will give a generalization of the convex hull
property for minimal hypersurfaces in hyperbolic space to CMC hypersurfaces. Similar results have been
obtained by Tonegawa in [To] and Alencar-Rosenberg in [AR]. We define a new notion $H$-shifted convex hull
which is obtained by shifting the convex hull up or down in a suitable way.

\vspace{0.3cm}

\noindent \textbf{Theorem 3.2.} [H-shifted Convex Hull Property for CMC Hypersurfaces] Let $-1<H<1$, and
$\Sigma$ be a CMC hypersurface with mean curvature $H$ in $\Hyp$ where $\PI\Sigma = \Gamma$. Then $\Sigma$ is
in the $H$-shifted convex hull of $\Gamma$.

\vspace{0.3cm}

We naturally generalize the definition of area minimizing hypersurfaces to CMC hypersurfaces, and call them
as "minimizing CMC hypersurfaces". These hypersurfaces are CMC hypersurfaces, and minimize area with the
volume constraint. Then, we adapt the exchange roundoff trick of Meeks-Yau for area minimizing hypersurfaces
to minimizing CMC hypersurfaces.

\vspace{0.3cm}

\noindent \textbf{Theorem 3.3.}[Exchange Roundoff Trick for Minimizing CMC Hypersurfaces] Let $\Gamma_1$ and
$\Gamma_2$ be two disjoint codimension-$1$ closed manifolds in $\SI$. If $\Sigma_1$ and $\Sigma_2$ are
minimizing CMC hypersurfaces with mean curvature $H$ in $\Hyp$ where $\PI \Sigma_i = \Gamma_i$, then
$\Sigma_1$ and $\Sigma_2$ are disjoint, too.

\vspace{0.3cm}

On the other hand, we will give a generic uniqueness result for minimizing CMC hypersurfaces in $\Hyp$.

\vspace{0.3cm}

\noindent \textbf{Corollary 4.3.} [Generic Uniqueness of Minimizing CMC hypersurfaces] Let $A$ be the space
of codimension-$1$ closed submanifolds of $\SI$, and let $A'\subset A$ be the subspace containing the closed
submanifolds of $\SI$ bounding a unique minimizing CMC hypersurface with mean curvature $H$ in $\Hyp$. Then
$A'$ is generic in $A$, i.e. $A-A'$ is a set of first category.

\vspace{0.3cm}

The organization of the paper as follows. In Section 2, we will give the basic definitions and results which
we use throughout the paper. In Section 3, we will generalize the convex hull property and exchange roundoff
trick to CMC hypersurfaces context. In Section 4, by using the results of Section 3, we will prove the
generic uniqueness result for minimizing CMC hypersurfaces. Finally, in Section 5, we will have some
concluding remarks.

\section{Preliminaries}

In this section, we will overview the basic results which we use in the following sections. Let $\Sigma^n$ be
a compact hypersurface, bounding a domain $\Omega^{n+1}$ in some ambient Riemannian manifold. Let $A$ be the
area of $\Sigma$, and $V$ be the volume of $\Omega$. Let's vary $\Sigma$ through a one parameter family
$\Sigma_t$, with corresponding area $A(t)$ and volume $V(t)$. If $f$ is the normal component of the
variation, and $H$ is the mean curvature of $\Sigma$, then we get $A'(0) = -\int_\Sigma n H f$, and
$V'(0)=\int_\Sigma f$ where $n$ is the dimension of $\Sigma$, and $H$ is the mean curvature.

Now let's define a new functional as a combination of $A$ and $V$. Let $I_H(t)= A(t) + n H V(t)$. Note that
$I_0(t)=A(t)$. If $\Sigma$ is a critical point of the functional $I_H$ for any variation $f$, then this will
imply $\Sigma$ has constant mean curvature $H$ [Gu].

If $\Sigma$ is hypersurface with boundary $\Gamma$, then we fix a hypersurface $M$ with $\partial M =
\Gamma$, and define $V(t)$ to be the volume of the domain bounded by $M$ and varied hypersurface. Again, if
$\Sigma$ is a critical point of the functional $I_H$ for any variation $f$, then this will imply $\Sigma$ has
constant mean curvature $H$. Note that critical point of the functional $I_H$ is independent of the choice of
the hypersurface $M$ since if $\widehat{I}_H$ is the functional which is defined with a different
hypersurface $\widehat{M}$, then $I_H - \widehat{I}_H = C$ for some constant $C$. In particular, $H=0$ is the
special case called minimal case, and the theory is very well-developed for this case [Ni], [CM].

\begin{defn} $[$ Minimal Case $]$
$\Sigma$ is called as \textit{minimal hypersurface} if it is critical point of $I_0$ (Area Functional)
for any variation. Equivalently, $\Sigma$ has constant mean curvature $0$ at every point. If also
$I_0''(0)>0$ for $\Sigma$ for any variation, then $\Sigma$ is called as \textit{stable minimal hypersurface}.
In this case, $\Sigma$ is locally area minimizing in the sense that $\Sigma$ has a small neighborhood $N$ in
the ambient manifold, and it has the smallest area among hypersurfaces in $N$ with same boundary. Finally,
$\Sigma$ is called \textit{area minimizing hypersurface} if $\Sigma$ is the absolute minimum of the
functional $I_0$ (having the smallest area) among hypersurfaces with same boundary. Clearly, all area
minimizing hypersurfaces are stable.
\end{defn}

For general $H$, the theory is called the constant mean curvature (CMC) case, and it is developed following
the traces of minimal case.

\begin{defn} $[$ CMC case $]$
$\Sigma$ is called as \textit{CMC hypersurface} if it is critical point of $I_H$ for any variation.
Equivalently, $\Sigma$ has constant mean curvature $H$ at every point. If also $I_H''(0)>0$ for $\Sigma$ for
any variation, then $\Sigma$ is called as \textit{stable CMC hypersurface}. In this case, $\Sigma$ is also
locally area minimizing in the sense that $\Sigma$ has a small neighborhood $N$ in the ambient manifold, and
it has the smallest area with a volume constraint among hypersurfaces in $N$ with same boundary [Br].
Following the minimal case, we will call $\Sigma$ as \textit{minimizing CMC hypersurface} if $\Sigma$ is the
absolute minimum of the functional $I_H$ among hypersurfaces with same boundary. Again, all minimizing CMC
hypersurfaces are stable.
\end{defn}

\noindent \textbf{Notation: }From now on, we will call CMC hypersurfaces with mean curvature $H$ as
\textit{$H$-hypersurfaces}.

\vspace{0.3cm}

We will call any noncompact hypersurface as stable (minimizing) $H$-hypersurface if any compact
codimension-$0$ submanifold with boundary is stable (minimizing) $H$-hypersurface.

Note that \textit{minimizing} $H$-hypersurfaces are natural generalizations of area minimizing hypersurfaces.
Even though this objects are widely used in the literature, they were not called with a special name before.
So, we will call them as minimizing $H$-hypersurfaces for our purposes. In the following sections, we will
also show that they have similar features with the area minimizing hypersurfaces.

After these general definitions of minimal and $H$-hypersurfaces, we will quote some basic facts about the
$H$-hypersurfaces in hyperbolic space. In this paper, we are interested in the complete $H$-hypersurfaces
$\Sigma^n$ in $\BH^{n+1}$ asymptotic to a closed codimension-$1$ submanifold $\Gamma^{n-1}$ in $\SI$. Now, we
will overview some basic results about these hypersurfaces.

First, we fix $\Gamma$ as a codimension-$1$ closed submanifold in $\SI$. $\Gamma$ separates $\SI$ into two
parts, say $\Omega^+$ and $\Omega^-$. By using these domains, we will give orientation to hypersurfaces in
$\Hyp$ asymptotic to $\Gamma$. With this orientation, mean curvature $H$ is positive if the mean curvature
vector points towards positive side of the hypersurface, negative otherwise. The following fact is a known as
maximum principle.

\begin{lem} $[$ Maximum Principle $]$ Let $\Sigma_1$ and $\Sigma_2$ be two hypersurfaces in a Riemannian manifold,
and intersect at a common point tangentially. If $\Sigma_2$ lies in positive side of $\Sigma_1$ around the
common point, then $H_1$ is strictly less than $H_2$ ($H_1 < H_2$) where $H_i$ is the mean curvature of
$\Sigma_i$ at the common point.
\end{lem}

With a simple application of this maximum principle, one can get the following well-known fact about
$H$-hypersurfaces in $\Hyp$ asymptotic to $\Gamma$.

\begin{lem}
Let $\Sigma$ be a complete $H$-hypersurface in $\Hyp$ asymptotic to a codimension-$1$ submanifold $\Gamma$ of
$\SI$. Then $|H|<1$.
\end{lem}

\begin{pf}
Let $\Sigma$ be as stated. Assume $H$ is positive. Let $x$ be a any point in $\Omega^+\subset\SI$. Let
$S_x(t)$ be the horospheres tangent to $\SI$ at $x$ parametrized so that $S_x(t)\rightarrow x$ as
$x\rightarrow 0$. Let $t$ increase until $S_x(t)$ touches $\Sigma$. Since $S_x(t)$ has mean curvature $1$, by
maximum principle $H<1$. If $H$ is negative, we can use same idea by choosing a point from $\Omega^-$. This
yields $H> -1$, and the proof follows.
\end{pf}

The following existence theorem for minimizing $H$-hypersurfaces in $\Hyp$ asymptotic to $\Gamma$ for a given
codimension-$1$ closed submanifold in $\SI$ was proved by Tonegawa [To], and Alencar-Rosenberg [AR]
independently by using geometric measure theory methods.

\begin{thm} $[$ Existence and Regularity $]$
Let $\Gamma$ be a codimension-$1$ closed submanifold in $\SI$, and let $|H|<1$. Then there exist a minimizing
$H$-hypersurface $\Sigma$ in $\Hyp$ where $\PI \Sigma = \Gamma$. Moreover, any such minimizing
$H$-hypersurface is smooth except a closed singularity set of dimension at most $n-7$.
\end{thm}

\section{Generalized Convex Hull property and Exchange Roundoff Trick}

In this section, we will give a generalization of the convex hull property of minimal hypersurfaces in $\Hyp$
to the $H$-hypersurfaces. Then, we will show a generalization of exchange roundoff trick for area minimizing
hypersurfaces to minimizing $H$-hypersurfaces.

Now, we define the convex hull of a subset $A$ of $\SI$ in $\Hyp$. If $\gamma$ is a round $n-1$-sphere in
$\SI$, then there is a unique geodesic plane $P$ in $\Hyp$ asymptotic to $\gamma$. $\gamma$ separates $\SI$
into two parts $\Omega^+$ and $\Omega^-$. Similarly, $P$ divides $\Hyp$ into two halfspaces $D^+$ and $D^-$
with $\PI D^\pm = \Omega^\pm$. We will call the halfspace whose asymptotic boundary containing $A$ as
\textit{supporting halfspace}. i.e. if $A\subset\Omega^+$, then $D^+$ is a supporting halfspace.

\begin{defn} $[$ Convex Hull $]$ Let $A$ be a subset of $\SI$. Then the \textit{convex hull} of $A$, $CH(A)$, is the smallest closed convex
subset of $\Hyp$ which is asymptotic to $A$. Equivalently, $CH(A)$ can be defined as the intersection of all
supporting closed half-spaces of $\Hyp$ [EM].
\end{defn}

Note that the asymptotic boundary of the convex hull of a subset of $\SI$ is the subset itself, i.e.
$\PI(CH(A)) = A$. In general, we say the hypersurface $\Sigma$ has the convex hull property if it is in the
convex hull of its boundary, i.e. $\Sigma \subset CH(\partial\Sigma)$. In special case, if $\Sigma$ is a
complete and noncompact hypersurface in $\Hyp$, then we say $\Sigma$ has convex hull property if it is in the
convex hull of its asymptotic boundary, i.e. $\Sigma\subset CH(\PI\Sigma)$. The following lemma is
well-known. The minimal hypersurfaces in $\Hyp$ have convex hull property [A1].

\begin{lem}$[$ Convex Hull Property for Minimal Hypersurfaces $]$
Let $\Sigma$ be a minimal hypersurface in $\Hyp$ with $\PI\Sigma = \Gamma$. Then $\Sigma \subset CH(\Gamma)$.
\end{lem}

\begin{pf} Let $\Sigma$ be a minimal hypersurface in $\Hyp$ with $\PI\Sigma = \Gamma$. By definition, if we
show that $\Sigma$ is disjoint from all (nonsupporting) halfspaces whose asymptotic boundary is disjoint from
$\Gamma$, then we are done.

Let $K$ be a nonsupporting halfspace in $\Hyp$, i.e. $\PI K \cap \Gamma = \emptyset$. Since $K$ is halfspace,
we can foliate $K$ with geodesic planes $P_t$ whose asymptotic boundaries are in $\PI K$. If $\Sigma$
intersect $K$, then $\Sigma$ must intersect a geodesic plane $P_{t_0}$ tangentially, as $\PI K \cap \Gamma =
\emptyset$. But since $\Sigma$ and $P_{t_0}$ are minimal, this contradicts to the maximum principle.
\end{pf}

Now, we will generalize this result to $H$-hypersurfaces in $\Hyp$. Of course, when $H\neq 0$,
$H$-hypersurfaces cannot have convex hull property (consider a $H$-hypersurface with asymptotic boundary is a
round sphere). But, we will modify this hypothesis naturally to suit for $H$-hypersurfaces.

If $\gamma$ is a round $n-1$-sphere in $\SI$, the geodesic plane $P$ asymptotic to $\gamma$ would be the
unique minimal hypersurface asymptotic to $\gamma$. Similarly, $R$-equidistant planes $\Sigma_{\pm R}$ to the
geodesic plane $\Sigma$ would be the unique $\pm H_R$-hypersurfaces asymptotic $\gamma$ [NS].  While $R$
varies from $0$ to $\infty$, $H_R$ varies from $0$ to $1$. We will call the $R$-equidistant plane to the
geodesic plane as $\Sigma_R$, if it is in the positive side of the geodesic plane, and as $\Sigma_{-R}$ if it
is in the negative side of the geodesic plane. By following the notation of the Section 2, clearly the mean
curvature of $\Sigma_R$ is $-H_R$ everywhere, and the mean curvature of $\Sigma_{-R}$ is $H_R$ everywhere.

Now, we modify the convex hull definition to generalize a similar notion for $H$-hypersurfaces. Let $\Gamma$
be a codimension-$1$ submanifold of $\SI$. Then, $\Gamma$ separates $\SI$ into two parts, say
$\Omega^\pm_\Gamma$ by giving an orientation. Now, we want orient all round $n-1$-spheres disjoint from
$\Gamma$ compatible with the orientation of $\Gamma$. Let $T$ be a round $n-1$-sphere in $\SI$ disjoint from
$\Gamma$. $T$ also separates $\SI$ into two parts. If $T\subset \Omega^+_\Gamma$, then call the part belong
in $\Omega^+_\Gamma$ as $\Omega^+_T$. If $T\subset \Omega^-_\Gamma$, then call the part belong in
$\Omega^-_\Gamma$ as $\Omega^-_T$. Clearly, this gives a compatible orientation to all such spheres.

Now, fix $\Gamma$ and orient all spheres accordingly. If $T$ is a round $n-1$-sphere in $\SI$, then there is
a unique $H$-hypersurface $P_H$ in $\Hyp$ asymptotic to $T$ for $-1<H<1$. $T$ separates $\SI$ into two parts
$\Omega^+$ and $\Omega^-$. Similarly, $P_H$ divides $\Hyp$ into two domains $D_H^+$ and $D_H^-$ with $\PI
D_H^\pm = \Omega^\pm$. We will call these regions as \textit{$H$-shifted halfspaces}. If the asymptotic
boundary $H$-shifted halfspace contains $\Gamma$, then we will call this $H$-shifted halfspace as
\textit{supporting $H$-shifted halfspace}. i.e. if $A\subset\Omega^+$, then $D_H^+$ is a supporting
$H$-shifted halfspace.

\begin{defn} $[$ H-shifted Convex Hull $]$ Let $\Gamma$ be a codimension-$1$ submanifold of $\SI$.
Then the $H$-shifted convex hull of $\Gamma$, $CH_H(\Gamma)$ is defined as the intersection of all supporting
closed $H$-shifted halfspaces of $\Hyp$.
\end{defn}

\begin{rmk}
Here, because of the orientation issue, we can define the $H$-shifted convex hull for only codimension-$1$
submanifolds of $\Hyp$. On the other hand, $0$-shifted convex hull is clearly same as the usual convex hull.
Intuitively, for $H>0$, $H$-shifted convex hull is obtained by shifting the convex hull into negative side,
and for $H<0$, it is obtained by shifting the convex hull into positive side. See Figure 1.
\end{rmk}

\begin{figure}[t]
\mbox{\vbox{\epsfbox{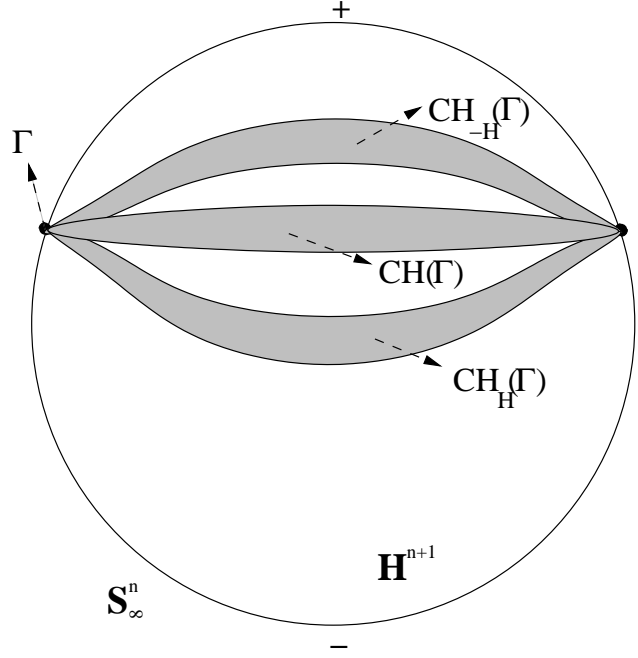}}} \caption{\label{fig:figure1} {$\Gamma$ is a codimension-$1$ closed
submanifold of $\SI$. $CH(\Gamma)$ is the convex hull of $\Gamma$. $CH_H(\Gamma)$ is the $H$-shifted convex
hull, and $CH_{-H}(\Gamma)$ is the $-H$-shifted convex hull of $\Gamma$ with the given orientation.}}
\end{figure}

Now, we will give a generalization of convex hull property of minimal hypersurfaces in $\Hyp$ to
$H$-hypersurfaces in $\Hyp$.

\begin{thm} $[$ H-shifted convex hull property for $H$-hypersurfaces $]$
Let $-1<H<1$, and $\Sigma$ be a $H$-hypersurface in $\Hyp$ where $\PI\Sigma = \Gamma$. Then $\Sigma$ is in
the $H$-shifted convex hull of $\Gamma$, i.e. $\Sigma \subset CH_H(\Gamma)$.
\end{thm}

\begin{pf} Let $\Sigma$ be a $H$-hypersurface in $\Hyp$ where $\PI\Sigma = \Gamma$.
By definition, if we show that $\Sigma$ is disjoint from all (nonsupporting) $H$-shifted halfspaces whose
asymptotic boundary is disjoint from $\Gamma$, then we are done.

Let $K$ be a $H$-shifted halfspace in $\Hyp$ with $\PI K \cap \Gamma = \emptyset$. Since $K$ is $H$-shifted
halfspace, we can foliate $K$ with $H$-hypersurfaces $P_t$ where their asymptotic boundaries are in $\PI K$.
If $\Sigma$ intersect $K$, then $\Sigma$ must intersect a $H$-hypersurfaces $P_{t_0}$ tangentially, as $\PI K
\cap \Gamma = \emptyset$. But since $\Sigma$ and $P_{t_0}$ have mean curvature $H$ everywhere, this
contradicts to the maximum principle (Lemma 2.1).
\end{pf}

\begin{rmk} This theorem is a natural generalization of the convex hull property for minimal hypersurfaces in
hyperbolic space to $H$-hypersurfaces in hyperbolic space. Similar versions of this result have been proved
by Alencar-Rosenberg in [AC], and by Tonegawa in [To].
\end{rmk}

Now, we are going to generalize exchange roundoff trick of Meeks-Yau for area minimizing hypersurfaces to
$H$-hypersurfaces. This technique mainly says that two area minimizing hypersurfaces with disjoint boundaries
are disjoint [Co3]. We will generalize this result to $H$-hypersurfaces in $\Hyp$.

\begin{thm}
Let $\Gamma_1$ and $\Gamma_2$ be two disjoint codimension-$1$ closed manifolds in $\SI$. If $\Sigma_1$ and
$\Sigma_2$ are minimizing $H$-hypersurfaces in $\Hyp$ where $\PI \Sigma_i = \Gamma_i$, then $\Sigma_1$ and
$\Sigma_2$ are disjoint, too.
\end{thm}

\begin{pf}
Assume that the minimizing $H$-hypersurfaces are not disjoint, i.e. $\Sigma_1\cap\Sigma_2\neq\emptyset$.
Since $\Sigma_i$ is a codimension-$1$ submanifold in $\Hyp$ whose asymptotic boundary is codimension-$1$
closed submanifold of $\SI$, $\Sigma_i$ separates $\Hyp$ into two parts. So, say $\Hyp-\Sigma_i =
\Omega^+_i\cup\Omega^-_i$.

Now, consider the intersection of hypersurfaces $\alpha=\Sigma_1\cap\Sigma_2$. Since the asymptotic
boundaries $\Gamma_1$ and $\Gamma_2$ are disjoint in $\SI$, by using the regularity results of Tonegawa in
[To], we can conclude that the intersection set $\alpha$ is in compact part of $\Hyp$. Moreover, by maximum
principle, the intersection cannot have isolated tangential intersections.

Now, without loss of generality, we assume that $\Sigma_1$ is {\em above} $\Sigma_2$ (the noncompact part of
$\Sigma_1$ lies in $\Omega^+_2$). Now define the compact hypersurfaces $S_i$ in $\Sigma_i$ as
$S_1=\Sigma_1\cap \Omega^-_2$, and $S_2 = \Sigma_2\cap\Omega^+_1$.  In other words, $S_1$ is the part of
$\Sigma_1$ lying {\em below} $\Sigma_2$, and $S_2$ is the part of $\Sigma_2$ lying {\em above} $\Sigma_1$.
Then, $\partial S_1 =\partial S_2 =\alpha$.

On the other hand, since $\Sigma_1$ and $\Sigma_2$ are minimizing $H$-hypersurfaces, then by definition, so
are $S_1$ and $S_2$, too. Then by swaping the hypersurfaces, we get new hypersurfaces $\Sigma_1 '$ and
$\Sigma_2 '$. In other words, let $\Sigma_1 ' = \{\Sigma_1-S_1\}\cup S_2$, and $\Sigma_2 ' =
\{\Sigma_2-S_2\}\cup S_1$ are new hypersurfaces with a singular set $\alpha$. See Figure 2.

We claim that $\Sigma_1 '$ and $\Sigma_2 '$ are also minimizing $H$-hypersurfaces. By Definition 2.2, we only
need to check that any compact part of $\Sigma_i '$ is minimum for the functional $I_H$ among the compact
hypersurfaces with same boundary. But since both $S_1$ and $S_2$ are minimizing $H$-hypersurfaces with same
boundary $\alpha$, then $I_H(S_1)=I_H(S_2)$. Let $T$ be a compact codimension-$0$ submanifold of $\Sigma_1$,
and $T'$ be the corresponding compact codimension-$0$ submanifold of $\Sigma_1 '$ with same boundary, i.e.
$T'=\{T-S_1\}\cup S_2$. Since $I_H(S_1)=I_H(S_2)$, this would imply $I_H(T)=I_H(T')$. Since $T$ is
minimizing, so is $T'$. This shows $\Sigma_1 '$ is minimizing $H$-hypersurface, and similarly $\Sigma_2 '$ is
also minimizing $H$-hypersurface.

\begin{figure}[t]
\mbox{\vbox{\epsfbox{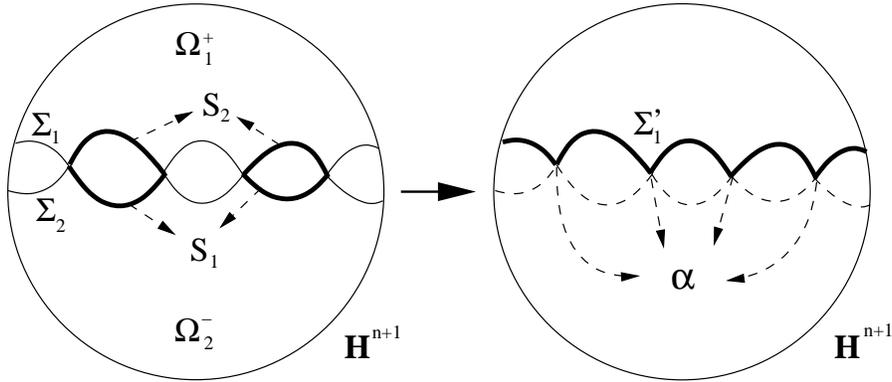}}} \caption{\label{fig:figure2} {$S_1$ is the part of $\Sigma_1$ lying below
$\Sigma_2$, and $S_2$ is the part of $\Sigma_2$ lying above $\Sigma_1$. After swaping $S_1$ and $S_2$, we get
a new minimizing $H$-hypersurface $\Sigma_1 '$ with singularity along $\alpha=\Sigma_1\cap\Sigma_2$.}}
\end{figure}

Now, in this new hypersurfaces $\Sigma_i '$, we have a codimension-$2$ singularity set along $\alpha$, which
contradicts to the regularity theorem for minimizing $H$-hypersurfaces, i.e. Theorem 2.3 (More intuitive way
to see this fact with Meeks-Yau's technique, one can modify the $\Sigma_1 '$ along the transverse
intersection on $\alpha$ in a suitable way, and one can show that this move reduces the functional $I_H$.).
\end{pf}

\section{Generic Uniqueness of $H$-hypersurfaces in $\Hyp$}

In this section, we will show that a generic codimension-$1$ closed submanifold in $\SI$ bounds a unique
$H$-hypersurface in $\Hyp$. The main idea is to adapt the technique in [Co3] to the $H$-hypersurfaces
context.

\begin{lem}
Let $\Gamma$ be a codimension-$1$ closed submanifold of $\SI$. Then either there exists a unique minimizing
$H$-hypersurface $\Sigma$ in $\Hyp$ asymptotic to $\Gamma$, or there are two canonical disjoint extremal
minimizing $H$-hypersurfaces $\Sigma^+$ and $\Sigma^-$ in $\Hyp$ asymptotic to $\Gamma$.
\end{lem}

\begin{pf}
Let $\Gamma$ be a closed submanifold of $\SI$. Since $\Gamma$ is a closed codimension-$1$ submanifold of the
sphere $\SI$, $\Gamma$ separates $\SI$ into two parts, say $\Omega^+$ and $\Omega^-$. Define sequences of
pairwise disjoint closed submanifolds of same topological type $\{\Gamma_i^+\}$ and $\{\Gamma_i^-\}$ in $\SI$
such that $\Gamma_i^+\subset \Omega^+$, and $\Gamma_i^-\subset \Omega^-$ for any $i$, and $\Gamma_i^+
\rightarrow \Gamma$, and $\Gamma_i^- \rightarrow \Gamma$ in Hausdorff metric. In other words,
$\{\Gamma_i^+\}$ and $\{\Gamma_i^-\}$ converges to $\Gamma$ from opposite sides.

By Theorem 2.3, for any $\Gamma_i^+\subset \SI$, there exist a minimizing $H$-hypersurface $\Sigma_i^+$ in
$\Hyp$. This defines a sequence of minimizing $H$-hypersurfaces $\{\Sigma_i^+\}$. Now, by using the sequence
$\{\Sigma_i^+\}$, define a new sequence of compact minimizing $H$-hypersurfaces $S_i^+\subset\Sigma_i^+$ with
$\partial S_i^+\rightarrow \Gamma\subset\SI$ (For example, let $S_i^+ =\Sigma_i^+\cap N_k(CH_H(\Gamma))$,
intersection of each $\Sigma_i^+$ with the $k$-neighborhood of the $H$-shifted convex hull of $\Gamma$ in
$\Hyp$ for sufficiently large $k$.) Then by using the proof of Theorem 2.3 (compactness theorem from
geometric measure theory), for this new sequence, we get a convergent subsequence $S_{i_j}^+\rightarrow
\Sigma^+$, and we get the minimizing $H$-hypersurface $\Sigma^+$ in $\Hyp$ asymptotic to $\Gamma$. Similarly,
we get the minimizing $H$-hypersurface $\Sigma^-$ in $\Hyp$ asymptotic to $\Gamma$. Now, we claim that these
minimizing $H$-hypersurfaces $\Sigma^\pm$ are very special by their construction, and they are disjoint from
each other.

Assume that $\Sigma^+$ and $\Sigma^-$ are not disjoint. Since these are minimizing $H$-hypersurfaces,
nontrivial intersection implies some part of $\Sigma^-$ lies {\em above} $\Sigma^+$. Since $\Sigma^+=\lim
S_i^+$, $\Sigma^-$ must also intersect some $S_i^+$ for sufficiently large $i$. But $S_i^+\subset
\Sigma_i^+$, and by Theorem 3.3, $\Sigma_i^+$ is disjoint from $\Sigma^-$ as $\PI\Sigma_i^+ = \Gamma_i^+$ is
disjoint from $\PI \Sigma^- = \Gamma$. This is a contradiction. This shows $\Sigma^+$ and $\Sigma^-$ are
disjoint.

Similar arguments show that $\Sigma^\pm$ are disjoint from any minimizing $H$-hypersurface $\Sigma '$
asymptotic $\Gamma$. Moreover, same argument shows that any minimizing $H$-hypersurface asymptotic to
$\Gamma$ must belong to the region bounded by $\Sigma^+$ and $\Sigma^-$ in $\Hyp$.
\end{pf}

\begin{rmk}
By above theorem and its proof, if $\Gamma$ bounds more than one minimizing $H$-hypersurface, then there
exist a canonical region $N$ in $\Hyp$ asymptotic to $\Gamma$ such that $N$ is the region between the
canonical minimizing $H$-hypersurfaces $\Sigma^+$ and $\Sigma^-$ in $\Hyp$. Moreover, by using similar ideas,
one can show that any minimizing $H$-hypersurface in $\Hyp$ asymptotic to $\Gamma$ is in the region $N$.
\end{rmk}

Now, we can prove the generic uniqueness result.

\begin{thm}
Let $B$ be the space of codimension-$1$ closed submanifolds of $\SI$ of same topological type, and let
$B'\subset B$ be the subspace containing the closed submanifolds of $\SI$ bounding a unique minimizing
$H$-hypersurface in $\Hyp$. Then $B'$ is generic in $B$, i.e. $B-B'$ is a set of first category.
\end{thm}

\begin{pf}
Fix a closed $n-1$-dimensional manifold $M$. Let $B= \{\Gamma\in C^0(M,S^n)\ | \ \Gamma(M) \mbox{ is an
embedding}\}$. Clearly, $B$ is an open subspace of $C^0(M,S^n)$. We will prove the theorem in 2 steps.

\vspace{0.3cm}

\textbf{Claim 1:} $B'$ is dense in $B$ as a subspace of $C^0(M,S^n)$ with the supremum metric.

\begin{pf}Let $B$ be the space of codimension-$1$ closed submanifolds of $\SI$ as described above. Let $\Gamma_0\in B$
be a closed submanifold in $\SI$. Since $\Gamma_0$ is closed submanifold, there exist a small regular
neighborhood $N(\Gamma_0)$ of $\Gamma_0$ in $\SI$, which is homeomorphic to $M\times I$. Let
$\Gamma:(-\epsilon,\epsilon)\rightarrow B$ be a small path in $B$ through $\Gamma_0$ such that
$\Gamma(t)=\Gamma_t$ and $\{\Gamma_t\}$ foliates $N(\Gamma_0)$ with closed submanifolds homeomorphic to $M$.
In other words, $\{\Gamma_t\}$ are pairwise disjoint closed submanifolds homeomorphic to $M$, and
$N(\Gamma_0)=\bigcup_{t\in (-\epsilon,\epsilon)} \Gamma_t$.

Since $\Gamma_0$ is a closed codimension-$1$ submanifold of $\SI$, $N(\Gamma_0)$ separates $\SI$ into two
parts, say $\Omega^+$ and $\Omega^-$, i.e. $\SI=N(\Gamma_0)\cup \Omega^+\cup \Omega^-$. Let $p^+$ be a point
in $\Omega^+$ and let $p^-$ be a point in $\Omega^-$ such that for a small $\delta$, $B_\delta(p^\pm)$ are in
the interior of $\Omega^\pm$. Let $\beta$ be the geodesic in $\Hyp$ asymptotic to $p^+$ and $p^-$.

By Lemma 4.1, for any $\Gamma_t$ either there exist a unique minimizing $H$-hypersurface $\Sigma_t$ in
$\Hyp$, or there is a canonical region $N_t$ in $\Hyp$ asymptotic to $\Gamma_t$, namely the region between
the canonical minimizing $H$-hypersurfaces $\Sigma_t^+$ and $\Sigma_t^-$. With abuse of notation, if
$\Gamma_t$ bounds a unique minimizing $H$-hypersurface $\Sigma_t$ in $\Hyp$, define $N_t=\Sigma_t$ as a
degenerate canonical neighborhood for $\Gamma_t$. Then, let $\widehat{N}= \{N_t\}$ be the collection of these
degenerate and nondegenerate canonical neighborhoods for $t\in(-\epsilon,\epsilon)$. Clearly, degenerate
neighborhood $N_t$ means $\Gamma_t$ bounds a unique minimizing $H$-hypersurface, and nondegenerate
neighborhood $N_s$ means that $\Gamma_s$ bounds more than one minimizing $H$-hypersurfaces. Note that by
Theorem 3.3, all canonical neighborhoods in the collection are pairwise disjoint. On the other hand, a
codimension-$1$ submanifold in $\Hyp$, whose asymptotic boundary is codimension-$1$ closed submanifold of
$\SI$, separates $\Hyp$ into two parts. So, the geodesic $\beta$ intersects all the canonical neighborhoods
in the collection $\widehat{N}$.

We claim that the part of $\beta$ which intersects $\widehat{N}$ is a finite line segment. Let $P^+$ be the
unique $H$-hypersurface asymptotic to the round sphere $\partial B_\delta(p^+)$ in $\Omega^+$. Similarly,
define $P^-$. By Theorem 3.3, $P^\pm$ are disjoint from the collection of canonical regions $\widehat{N}$.
Let $\beta\cap P^\pm=\{q^\pm\}$. Then the part of $\beta$ which intersects $\widehat{N}$ is the line segment
$l\subset \beta$ with endpoints $q^+$ and $q^-$. Let $C$ be the length of this line segment $l$. See Figure
3.

\begin{figure}[t]
\mbox{\vbox{\epsfbox{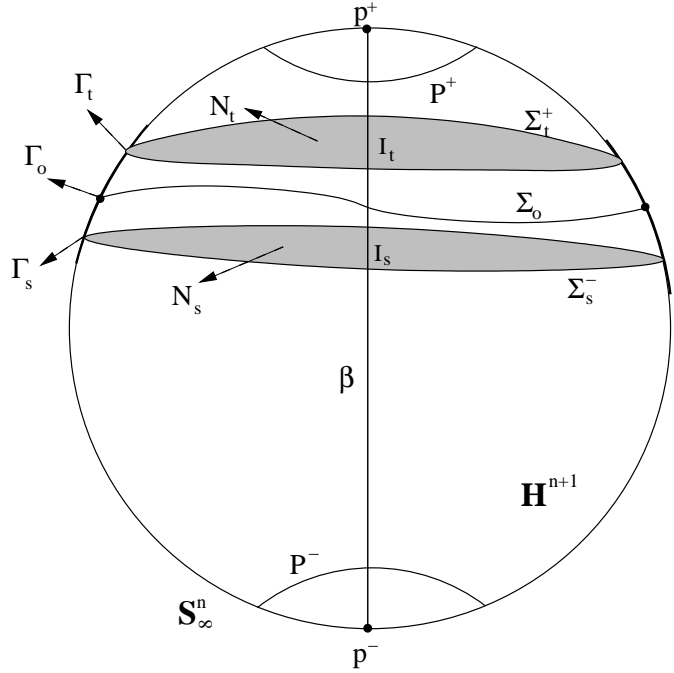}}} \caption{\label{fig:figure3} {A finite segment of geodesic $\gamma$
intersects the collection of minimizing $H$-hypersurfaces $\Sigma_t$ in $\BH^n$ asymptotic to $\Gamma_t$ in
$\SI$.}}
\end{figure}

Now, for each $t\in(-\epsilon,\epsilon)$, we will assign a real number $s_t\geq 0$. If there exists a unique
minimizing $H$-hypersurface $\Sigma_t$ for $\Gamma_t$ ($N_t$ is degenerate), then let $s_t$ be $0$. If not,
let $I_t = \beta\cap N_t$, and $s_t$ be the length of $I_t$. Clearly if $\Gamma_t$ bounds more than one
minimizing $H$-hypersurface ($N_t$ is nondegenerate), then $s_t > 0$. Also, it is clear that for any $t$,
$I_t\subset l$ and $I_t\cap I_s=\emptyset$ for any $t,s\in (-\epsilon,\epsilon)$. Then,
$\sum_{t\in(-\epsilon,\epsilon)} s_t < C$ where $C$ is the length of $l$. This means for only countably many
$t\in(-\epsilon,\epsilon)$, $s_t> 0$. So, there are only countably many nondegenerate $N_t$ for
$t\in(-\epsilon,\epsilon)$. Hence, for all other $t$, $N_t$ is degenerate. This means there exist uncountably
many $t\in(-\epsilon,\epsilon)$, where $\Gamma_t$ bounds a unique minimizing $H$-hypersurface. Since
$\Gamma_0$ is arbitrary, this proves $B'$ is dense in $B$.
\end{pf}

\textbf{Claim 2:} $B'$ is generic in $B$, i.e. $B-B'$ is a set of first category.

\begin{pf} We will prove that $B '$ is countable intersection of open dense subsets of a complete metric space. Then the
result will follow by Baire category theorem.

Since the space of continuous maps from $M$ to $n$-sphere $C^0(M,S^n)$ is complete with supremum metric, then
the closure of $B$ in $C^0(M,S^n)$, $\bar{B}\subset C^0(M,S^n)$, is also complete. Note that $B$ is an open
subspace of $C^0(M,S^n)$.

Now, we will define a sequence of open dense subsets $U^i\subset B$ such that their intersection will give us
$B '$. Let $\Gamma\in B$ be a closed codimension-$1$ submanifold of $\SI$ homeomorphic to $M$. Let
$N(\Gamma)\subset \SI$ be a regular neighborhood of $\Gamma$ in $\SI$, which is homeomorphic to $M \times I$.
Then, define an open neighborhood of $\Gamma$ in $B$, $U_\Gamma\subset B$, such that $U_\Gamma = \{\alpha \in
B \ | \ \alpha(M)\subset N(\Gamma)\}$. Clearly, $B= \bigcup_{\Gamma\in B} U_\Gamma$.  Now, define a geodesic
$\beta_\Gamma$ as in Claim 1, which intersects all the minimizing $H$-hypersurfaces asymptotic to
submanifolds in $U_\Gamma$.

Now, as in Claim 1, for any $\alpha \in U_\Gamma$, there exist a canonical region $N_\alpha$ in $\Hyp$ (which
can be degenerate if $\alpha$ bounds a unique $H$-hypersurface). Let $I_{\alpha,\Gamma} = N_\alpha \cap
\beta_\Gamma$. Then let $s_{\alpha,\Gamma}$ be the length of $I_{\alpha,\Gamma}$ ($s_{\alpha,\Gamma}$ is $0$
if $N_\alpha$ degenerate). Hence, for every element $\alpha$ in $U_\Gamma$, we assign a real number
$s_{\alpha,\Gamma} \geq 0$.

Now, we define the sequence of open dense subsets in $U_\Gamma$. Let $U^i_\Gamma = \{\alpha\in U_\Gamma \ | \
s_{\alpha,\Gamma} < 1/i \  \}$. We claim that $U^i_\Gamma$ is an open subset of $U_\Gamma$ and $B$. Let
$\alpha\in U^i_\Gamma$, and let $s_{\alpha,\Gamma} = \lambda < 1/i$. So, the interval
$I_{\alpha,\Gamma}\subset \beta_\Gamma$ has length $\lambda$. let $I ' \subset \beta_\Gamma$ be an interval
containing $I_{\alpha,\Gamma}$ in its interior, and has length less than $1/i$. By the proof of Claim 1, we
can find two submanifolds $\alpha^\pm\in U_\Gamma$ with the following properties.

\begin{itemize}

\item $\alpha^\pm$ are disjoint from $\alpha$,

\item $\alpha^\pm$ are lying in opposite sides of $\alpha$ in $\SI$,

\item $\alpha^\pm$ bounds unique minimizing $H$-hypersurface $\Sigma_{\alpha^\pm}$,

\item $\Sigma_{\alpha^\pm} \cap \beta_\Gamma \subset I '$.

\end{itemize}

The existence of such submanifolds is clear from the proof of Claim 1, as if one takes any foliation
$\{\alpha_t\}$ of a small neighborhood of $\alpha$ in $\SI$, there are uncountably many submanifolds in the
family bounding a unique minimizing $H$-hypersurface, and one can choose sufficiently close pair of
submanifolds to $\alpha$, to ensure the conditions above.

After finding $\alpha^\pm$, consider the open neighborhood homeomorphic to $M \times I$, $F_\alpha$, in $\SI$
bounded by $\alpha^+$ and $\alpha^-$. Let $V_\alpha = \{ \gamma\in U_\Gamma \ | \ \gamma(M)\subset F_\alpha
\}$. Clearly, $V_\alpha$ is an open subset of $U_\Gamma$. If we can show $V_\alpha\subset U^i_\Gamma$, then
this proves $U^i_\Gamma$ is open for any $i$ and any $\Gamma\in A$.

Let $\gamma\in V_\alpha$ be any submanifold, and $N_\gamma$ be its canonical neighborhood described in the
proof of Lemma 4.1 and in Remark 4.1. Since $\gamma(M)\subset F_\alpha$, $\alpha^+$ and $\alpha^-$ lie in
opposite sides of $\gamma$ in $\SI$. This means $\Sigma_{\alpha^+}$  and $\Sigma_{\alpha^-}$ lie in opposite
sides of $N_\alpha$. By choice of $\alpha^\pm$, this implies $N_\gamma \cap \beta_\Gamma= I_{\gamma,\Gamma}
\subset I '$. So, the length $s_{\gamma,\Gamma}$ is less than $1/i$. This implies $\gamma\in U^i_\Gamma$, and
so $V_\alpha\subset U^i_\Gamma$. Hence, $U^i_\Gamma$ is open in $U_\Gamma$ and $B$.

Now, we can define the sequence of open dense subsets. let $U^i = \bigcup_{\Gamma\in B} U^i_\Gamma$ be an
open subset of $B$. Since, the elements in $B '$ represent the submanifolds bounding a unique minimizing
$H$-hypersurface, for any $\alpha\in B '$, and for any $\Gamma\in B$, $s_{\alpha,\Gamma} = 0$. This means
$B'\subset U^i$ for any $i$. By Claim 1, $U^i$ is open dense in $B$ for any $i>0$.

As we mention at the beginning of the proof, since the space of continuous maps from $M$ to $n$-sphere
$C^0(M,S^n)$ is complete with supremum metric, then the closure of $B$ in $C^0(M,S^n)$, $\bar{B}$, is also
complete metric space. Since $B'$ is dense in $B$, it is also dense in $\bar{B}$. This implies $U^i$ is a
sequence of open dense subsets of $\bar{B}$. On the other hand, since $s_{\alpha,\Gamma} = 0$ for any
$\alpha\in B '$, and for any $\Gamma\in B$, $B '  = \bigcap_{i>0} U^i$. Then, $\bar{B}- B' = \bigcup_{i>0}
{U^i}^c$, where ${U^i}^c$ is complement of $U^i$. Since $U^i$ is open dense, then ${U^i}^c$ is nowhere dense.
Since $B-B'$ is countable union of nowhere dense sets, and $\bar{B}$ is a complete metric space, $B-B'$ is a
set of first category, by Baire Category Theorem. Hence, $B'$ is generic in $B$.
\end{pf}

\end{pf}

\begin{rmk}
This genericity result can be generalized to a genericity in all closed codimension-$1$ submanifolds of
$\SI$. We can stratify the whole space of codimension-1 closed submanifolds of $\SI$ by topological type, and
this result give us a genericity result in each strata. This implies genericity in the whole.
\end{rmk}

\begin{cor}
Let $A$ be the space of codimension-$1$ closed submanifolds of $\SI$, and let $A'\subset A$ be the subspace
containing the closed submanifolds of $\SI$ bounding a unique minimizing $H$-hypersurface in $\Hyp$. Then
$A'$ is generic in $A$, i.e. $A-A'$ is a set of first category.
\end{cor}

\section{Concluding Remarks}

In this paper, we showed a generic $n-1$-submanifold of $\SI$ bounds a unique \textit{minimizing}
$H$-hypersurface in $\Hyp$. The condition being minimizing in this result is necessary for our techniques.
This is because we need Theorem 3.3 for this method, and being just $H$-hypersurface is not enough for that
theorem. However, in the case $n=3$, one might get a similar result for general $H$-hypersurfaces by using
the methods in [Co1].

We consider CMC hypersurfaces as generalizations of minimal hypersurfaces. In most cases, the known results
on CMC hypersurfaces are some form of generalization of a result about minimal hypersurfaces. In other words,
the interesting questions about CMC hypersurfaces are mostly generalizations its correspondent in minimal
hypersurfaces world. In particular, one might look at the results in this paper as generalizations of some
notions about minimal hypersurfaces to CMC hypersurfaces in hyperbolic space.

On the other hand, in special case $H=0$, there is no known example of an asymptotic $n-1$-manifold in $\SI$
bounding more than one area minimizing hypersurface in $\Hyp$. So, when we generalize this question to CMC
hypersurfaces, we get the following problem: Is there any asymptotic $n-1$-manifold in $\SI$ bounding more
than one minimizing $H$-hypersurface in $\Hyp$? Or more boldly one can ask the following question:

\vspace{0.3cm}

\noindent \textbf{Question 1.} Is it true that for any codimension-$1$ closed submanifold $\Gamma^{n-1}$ of
$\SI$, there exist a unique minimizing $H$-hypersurface $\Sigma^n$ in $\Hyp$ asymptotic to $\Gamma$, or not?

\vspace{0.3cm}

Moreover, there are interesting problems when one think about the all minimizing $H$-hypersurfaces in $\Hyp$
asymptotic to same submanifold in $\SI$. Intuitively, one might think that two minimizing $H$-hypersurfaces
asymptotic to same submanifold in $\SI$ should be disjoint. In fact this is the case when they are compact
and they have same boundary. But, in noncompact case (even in the case $H=0$), this is not known.

\vspace{0.3cm}

\noindent \textbf{Question 2.} Let $\Gamma$ be a codimension-$1$ submanifold of $\SI$. If $\Sigma_1$ is a
minimizing $H_1$-hypersurface and $\Sigma_2$ is a minimizing $H_2$-hypersurface with $\PI \Sigma_i =\Gamma$,
then are $\Sigma_1$ and $\Sigma_2$ necessarily disjoint?

\vspace{0.3cm}

When we consider the all minimizing $H$-hypersurfaces in $\Hyp$ asymptotic to same submanifold $\Gamma$ in
$\SI$, the natural question one might ask the foliation of $\Hyp$ with such $H$-hypersurfaces. Chopp and
Velling raised the question "when do such $H$-hypersurfaces constitute a foliation of $\Hyp$?" in [CV] and
gave interesting examples for $\BH^3$ by using computational techniques. Note that existence of such a
foliation for $\Gamma$ implies uniqueness of $H$-hypersurfaces asymptotic to $\Gamma$ by maximum principle.
Since such a foliation implies uniqueness for each $H$, following the result of this paper, one might ask
whether these asymptotic submanifolds with such foliations are generic?

The answer to this question is "no". Such a foliation implies uniqueness of not only minimizing CMC
hypersurfaces, but also any type of CMC hypersurfaces. On the other hand, by Hass's result in [Ha] and
Anderson's result in [A2], there are examples of Jordan curves in $\Si$ bounding many minimal surfaces in
$\BH^3$. In particular, in Hass's construction, one can get a genus 1 minimal surfaces asymptotic to the
curve, and this implies existence of minimal planes in both sides of the surface by using Meeks-Yau's results
[MY]. But, since the curves which Hass constructed contains an open set of Jordan curves, this implies the
curves inducing foliation by $H$-hypersurfaces cannot be generic.

\end{document}